\documentclass[12pt,a4paper]{article}

\usepackage{amssymb}

\usepackage{amsmath}

\newtheorem{theorem}{Theorem}[section]
\newtheorem{definition}[theorem]{Definition}
\newtheorem{lemma}[theorem]{Lemma}

\begin{document}
\title{Rosso-Yamane Theorem on PBW basis of $U_q(A_N)$\footnote{Supported by the
NNSF of China (No.10771077) and the NSF of Guangdong Province
(No.06025062).}}
\author{Yuqun Chen, Hongshan Shao   \\
\\
{\small \ School of Mathematical Sciences}\\
{\small \ South China Normal University}\\
{\small \ Guangzhou 510631, P. R. China}\\
{\small \ yqchen@scnu.edu.cn}\\
{\small \ shaohongshan118@163.com} \\
\\
K. P. Shum\\
{\small \ Department of Mathematics}\\
{\small \  The University of Hong Kong}\\
{\small \ Pokfulam Road, Hong Kong, China (SAR)}\\
{\small \ kpshum@maths.hku.hk}\\
}
\date{}
\maketitle

\maketitle \noindent\textbf{Abstract:} Let $U_q(A_N)$ be the
Drinfeld-Jimbo quantum group of type $A_N$. In this paper, by using
Gr\"{o}bner-Shirshov bases, we give a simple (but not short) proof
of the Rosso-Yamane Theorem on PBW basis of  $U_q(A_N)$.

\maketitle \noindent\textbf{Key words: } Quantum group; Quantum
enveloping algebra; Gr\"{o}bner-Shirshov basis.

\noindent \textbf{AMS 2000 Subject Classification}: 20G42, 16S15,
13P10

\section{Introduction}

Since any algebra (commutative, associative, Lie), as well as any
module over an algebra, can be presented by generators and defining
relations, it is important to have a general method to deal with
these presentations. Such a method now exists and is called the
Gr\"{o}bner bases method (due to B. Buchberger \cite{Bu65},
\cite{Bu70}), or standand bases method (due to H. Hironaka
\cite{Hi}), or  Gr\"{o}bner-Shirshov bases method (due to A. I.
Shirshov \cite{Sh}). The original Shirshov's paper \cite{Sh} is for
Lie algebra presentations, but it can be easily adopted for
associative algebra presentations as well, see L. A. Bokut
\cite{B76} and G. Bergman \cite{Be}.

Let, for example, $L=Lie(X|[x_ix_j]-\sum\alpha_{ij}^kx_k, \ i>j, \
x_i,x_j,x_k \in X)$ be a Lie algebra over a field (or a commutative
ring) $k$ presented by a $k$-basis $X$ and the multiplication table.
Then $S=\{[x_ix_j]-\sum\alpha_{ij}^kx_k| \ i>j, \ x_i,x_j,x_k\in
X\}$ is a Gr\"{o}bner-Shirshov basis (subset) of the free Lie
algebra $Lie(X)$ over $k$. On the other hand, the universal
enveloping algebra $U(L)=k\langle
X|x_ix_j-x_jx_i-\sum\alpha_{ij}^kx_k, \ i>j, \ x_i,x_j,x_k \in
X\rangle$ is the associative algebra presented by the same set $X$
and the defining relations $S^{(-)}$ (we rewrite $S$ using
$[xy]=xy-yx$). There is a general but not difficult result that for
any $S\subset Lie(X), \ S$ is a Gr\"{o}bner-Shirshov basis in the
sense of Lie algebras if and only if $S^{(-)}\subset k\langle
X\rangle$ is a Gr\"{o}bner-Shirshov basis in  the sense of
associative algebras (see, for example, \cite{BKLM} and
\cite{BoCh1}). This means that in our case, $S^{(-)}$ is a
Gr\"{o}bner-Shirshov basis (subset) in $k\langle X\rangle$. By
Composition-Diamond lemma (see below), the $S$-irreducible words on
$X$, $Irr(S)=\{x_{i_1}\ldots x_{i_k}, \ i_1\leq\ldots \leq i_k, \
k\geq 0\}$ form a $k$-basis of $U(L)$. This is a conceptional proof
of the PBW-Theorem by using Gr\"{o}bner-Shirshov bases theory.

There are many results on Gr\"{o}bner-Shirshov bases for associative
and Lie algebras, as well as for semigroups, groups, conformal
algebras, dialgebras, and so on, see, for example, surveys
\cite{BoKo2}, \cite{BoKo3}, \cite{kl5} and \cite{BoCh2}. Let us
mention those for simple  Lie algebras and Lie superalgebras via
Serre's presentations (\cite{BoKl1}, \cite{BoKl2}, \cite{BoKl3},
\cite{BoKl4}, \cite{BKLM}), for modules over simple Lie algebras and
Iwahori-Hecke algebras (\cite{KaLe1}, \cite{KaLe2}, \cite{kl5}), for
Kac-Moody algebras of types $A_n^{(1)}, \ B_n^{(1)}, \ C_n^{(1)}, \
D_n^{(1)}$ (\cite{poro1}, \cite{poro2}, \cite{poro3}), for Coxeter
groups (\cite{BoSh}), for braid groups via Artin-Burau,
Artin-Garside and Briman-Ko-Lee presentations (\cite{B08}, \cite{Bo}
and \cite{BCS}).

Drinfeld-Jimbo (\cite{Dr}, \cite{Ji1}) presentations for quantized
enveloping algebras $U_q(g)$, where $g$ is a semisimple Lie algebra,
are a natural source  of associative presentations. M. Rosso
\cite{Ro} and I. Yamane \cite{Ya} found the PBW-basis of $U_q(A_N)$.
G. Lusztig \cite{Lu1} and \cite{Lu2}, and M. Kashiware \cite{Ka1}
and \cite{Ka2} found the bases of $U_q(g)$ for any semisimple
algebra $g$, as well as for their representations. Their approach
work equally well for quantized enveloping algebras associated with
arbitrary symmetrizable Cartan matrix, not just those corresponding
to finite dimensional Lie algebras.  V. K. Kharchenko \cite{Kh}
found the approach to linear bases of quantized enveloping algebras
via the so called character Hopf algebras.

It the paper \cite{BoMa},  Gr\"{o}bner-Shirshov bases approach was
applied to study $U_q(g)$ for any symmetrizable Cartan matrix. Using
this approach, they got a new proof of the triangular decomposition
of $U_q(g)$ (see, for example, Yantzen \cite{Ya2}). For $U_q(A_N)$,
it was proved by Bokut and Malcolmson \cite{BoMa} that the Jimbo
relations (see \cite{Ya}) of $U_q^+(A_N)$ constitute a
Gr\"{o}bner-Shirshov basis of $U_q^+(A_N)$ in Jimbo generators
$x_{ij}, 1\leq i,j\leq N+1$ (see below).

In this paper, we give an elementary proof that Jimbo relations $S$
is a Gr\"{o}bner-Shirshov basis of $U_q^+(A_N)$. For such a purpose,
we just check all possible compositions of polynomials from $S$ and
proved that all them can be resolved. Also in $\S 1$ in this paper,
we are giving necessary definitions and Composition-Diamond lemma
following Shirshov \cite{Sh}.

\section{Preliminaries}

We first cite some concepts and results from the literature which
are related to the Gr\"{o}bner-Shirshov bases for associative
algebras.

Let $k$ be a field, $k\langle X\rangle$ the free associative algebra
over $k$ generated by $X$ and $ X^{*}$ the free monoid generated by
$X$, where the empty word is the identity which is denoted by 1. For
a word $w\in X^*$, we denote the length of $w$ by $deg(w)$. Let
$X^*$ be a well ordered set. Let $f\in k\langle X\rangle$ with the
leading word $\bar{f}$. Then we call $f$  monic if $\bar{f}$ has
coefficient 1.

\begin{definition} (\cite{Sh}, see also \cite{B72}, \cite{B76}) \
Let $f$ and $g$ be two monic polynomials in \textmd{k}$\langle
X\rangle$ and $<$ a well order on $X^*$. Then, there are two kinds
of compositions:

$(i)$ If \ $w$ is a word such that $w=\bar{f}b=a\bar{g}$ for some
$a,b\in X^*$ with deg$(\bar{f})+$deg$(\bar{g})>$deg$(w)$, then the
polynomial
 $(f,g)_w=fb-ag$ is called the intersection composition of $f$ and
$g$ with respect to $w$.

$(ii)$ If  $w=\bar{f}=a\bar{g}b$ for some $a,b\in X^*$, then the
polynomial $(f,g)_w=f - agb$ is called the inclusion composition of
$f$ and $g$ with respect to $w$.

\end{definition}

\begin{definition}(\cite{B72}, \cite{B76}, cf. \cite{Sh})
Let $S\subset k\langle X\rangle$ such that every $s\in S$ is monic.
Then the composition $(f,g)_w$ is called trivial modulo $(S,w)$ if
$(f,g)_w=\sum\alpha_i a_i s_i b_i$, where each $\alpha_i\in k$,
$a_i,b_i\in X^{*}, \ s_i\in S$ and $\overline{a_i s_i b_i}<w$. If
this is the case, then we write
$$
(f,g)_w\equiv0\quad mod(S,w).
$$
In
general, for $p,q\in k\langle X\rangle$, we write $ p\equiv q\quad
mod(S,w) $ which means that $p-q=\sum\alpha_i a_i s_i b_i $, where
each $\alpha_i\in k,a_i,b_i\in X^{*}, \ s_i\in S$ and $\overline{a_i
s_i b_i}<w$.
\end{definition}

\begin{definition} (\cite{B72}, \cite{B76}, cf. \cite{Sh}) \
We call the set $S$ with respect to the well order $<$ a
Gr\"{o}bner-Shirshov set (basis) in $k\langle X\rangle$ if any
composition of polynomials in $S$ is trivial modulo $S$.
\end{definition}

If a subset $S$ of $k\langle X\rangle$ is not a Gr\"{o}bner-Shirshov
basis, then we can add to $S$ all nontrivial compositions of
polynomials of $S$, and by continuing this process (maybe
infinitely) many times, we eventually obtain a Gr\"{o}bner-Shirshov
basis $S^{comp}=S^c$. Such a process is called the Shirshov
algorithm.

A well order $>$ on $X^*$ is called monomial if it is compatible
with the multiplication of words, that is, for $u, v\in X^*$, we
have
$$
u > v \Rightarrow w_{1}uw_{2} > w_{1}vw_{2},  \ for \  all \
 w_{1}, \ w_{2}\in  X^*.
$$
A standard example of monomial order on $X^*$ is the deg-lex order
to compare two words first by degree and then lexicographically,
where $X$ is a linearly ordered set.

The following lemma was first proved by Shirshov \cite{Sh} for free
Lie algebras (with deg-lex order) in 1962 (see also Bokut
\cite{B72}). In 1976, Bokut \cite{B76} specialized the approach of
Shirshov to associative algebras (see also Bergman \cite{Be}). For
the case of commutative polynomials, this lemma is known as the
Buchberger's Theorem in \cite{Bu65} and \cite{Bu70}.

\begin{lemma}\label{l1}
(Composition-Diamond Lemma) \ Let $k$ be a field, $A=k \langle
X|S\rangle=k\langle X\rangle/Id(S)$ and $<$ a monomial order on
$X^*$, where $Id(S)$ is the ideal of $k \langle X\rangle$ generated
by $S$. Then the following statements are equivalent:
\begin{enumerate}
\item[(i)] $S $ is a Gr\"{o}bner-Shirshov basis.
\item[(ii)] $f\in Id(S)\Rightarrow \bar{f}=a\bar{s}b$
for some $s\in S$ and $a,b\in  X^*$.
\item[(iii)] $Irr(S) = \{ u \in X^* |  u \neq a\bar{s}b ,s\in S,a ,b \in X^*\}$
is a basis of the algebra $A=k\langle X | S \rangle. \ \square$
\end{enumerate}
\end{lemma}

\section{Rosso-Yamane theorem on PBW basis of $U_q(A_N)$}

Let $k$ be a field, $A=(a_{ij})$ an integral symmetrizable $N \times
N$ Cartan matrix so that $a_{ii}=2$, $a_{ij} \leq0 \ (i \neq j)$ and
there exists a diagonal matrix $D$ with diagonal entries $d_i$ which
are nonzero integers such that the product $DA$ is symmetric. Let q
be a nonzero element of $k$ such that $q^{4d_i} \neq 1$ for each
$i$. Then the quantum enveloping algebra is (see \cite{Dr},
\cite{Ji1})
$$
U_q(A)=k \langle X \cup H \cup Y|S^+ \cup K \cup T \cup S^- \rangle,
$$
where

\begin{eqnarray*}
X&=&\{x_i\},\\
H&=&\{h_i^{\pm1}\},\\
Y&=&\{y_i\},\\
S^+&=&\{ \sum_{\nu=0}^{1-a_{ij}}(-1)^\nu \left(
\begin{array}{c}
1-a_{ij} \\
\nu
\end{array}
 \right)_tx_i^{1-a_{ij}-\nu}x_jx_i^\nu, \ where \ i \neq j, \ t=q^{2d_i}\},\\
S^-&=& \{ \sum\limits_{\nu=0}^{1-a_{ij}}(-1)^\nu\left(
\begin{array}{c}
1-a_{ij} \\
\nu
\end{array}
\right)_ty_i^{1-a_{ij}-\nu}y_jy_i^\nu, \ where \ i \neq j, \ t=q^{2d_i}\},\\
K&=&\{h_ih_j-h_jh_i, h_ih_i^{-1}-1, h_i^{-1}h_i-1, x_jh_i^{\pm1}-q^{\mp1}
d_ia_{ij}h^{\pm1}x_j, h_i^{\pm1}y_j-q^{\mp1}y_jh^{\pm1}\},  \\
T&=&\{x_iy_j-y_jx_i- \delta_{ij}
\frac{h_i^2-h_i^{-2}}{q^{2d_i}-q^{-2d_i}}\} \ \ \mbox{ and }\\
\left(
\begin{array}{c}
m \\
 n \\
 \end{array}
 \right)_t&=&\left\{
 \begin{array}{cc}
 \prod\limits_{i=1}^n \frac{t^{m-i+1}-t^{i-m-1}}{t^i-t^{-i}} & \ \ (for \ m>n>0),\\
 1 & \ \ \ (for \ n=0 \ or \ m=n).
 \end{array}
 \right.
\end{eqnarray*}

Let
$$
A=A_N=\left(
  \begin{array}{ccccc}
    2 & -1 & 0 & \cdots & 0 \\
    -1 & 2 & -1 & \cdots & 0 \\
    0 & -1 & 2 & \cdots & 0 \\
    \cdot & \cdot & \cdot & \cdot & \cdot \\
    0 & 0 & 0 & \cdots & 2 \\
  \end{array}
\right) \ \mbox{ and }\ \ q^8 \neq1.
$$
It is reminded that in this case, the diagonal matrix $D$ is
identity.

We introduce some new variables defined by Jimbo (see \cite{Ya})
which generate $U_q(A_N)$:
$$
\widetilde{X}=\{x_{ij}, 1 \leq i<j \leq N+1\},
$$
where
$$
x_{ij}=\left\{
\begin{array}{cc}
x_i & \ \ \  j=i+1,\\
qx_{i,j-1}x_{j-1,j}-q^{-1}x_{j-1,j}x_{i,j-1} & \ \ \ j>i+1.
\end{array}
\right.
$$
We now order the set $\widetilde{X}$ in the following way.
$$
x_{mn}>x_{ij} \Longleftrightarrow(m,n)>_{lex}(i,j).
$$
Let us recall from Yamane \cite{Ya}  the following notation:
\begin{eqnarray*}
C_1&=&\{((i,j),(m,n))|i=m<j<n\},\\
C_2&=&\{((i,j),(m,n))|i<m<n<j\},\\
C_3&=&\{((i,j),(m,n))|i<m<j=n\},\\
C_4&=&\{((i,j),(m,n))|i<m<j<n\},\\
C_5&=&\{((i,j),(m,n))|i<j=m<n\},\\
C_6&=&\{((i,j),(m,n))|i<j<m<n\}.
\end{eqnarray*}
Let the set $\widetilde{S}^+$ consist of Jimbo relations:
\begin{eqnarray*}
x_{mn}x_{ij}&-&q^{-2}x_{ij}x_{mn}  \ \ \ \ \ \ \ \ \ \ \ \ \ \ \ \
\ \ \ \ \ \ \ \ \ ((i,j),(m,n)) \in C_1 \cup C_3,\\
x_{mn}x_{ij}&-&x_{ij}x_{mn} \ \ \ \ \ \ \ \ \ \ \ \ \ \ \ \ \ \ \
\ \ \ \ \ \ \ \ \ \ ((i,j),(m,n)) \in C_2
\cup C_6,\\
x_{mn}x_{ij}&-&x_{ij}x_{mn}+(q^2-q^{-2})x_{in}x_{mj} \ \ \
((i,j),(m,n)) \in C_4,\\
x_{mn}x_{ij}&-&q^2x_{ij}x_{mn}+qx_{in} \ \ \ \ \ \ \ \ \ \ \ \ \ \
\ \ \ ((i,j),(m,n)) \in C_5.
\end{eqnarray*}
It is easily seen that $U^+_q(A_N)=k\langle
\widetilde{X}|\widetilde{S^+}\rangle$.

The following theorem is from \cite{BoMa}.

\begin{theorem} (\cite{BoMa} Theorem 4.1)\label{t1}
Let the notation be as before. Then, with the deg-lex order on
$\widetilde{X}^*$, $\widetilde{S}^+$ is a Gr\"{o}bner-Shirshov basis
for $k\langle \widetilde{X}|\widetilde{S^+}\rangle=U^+_q(A_N)$.
\end{theorem}
{\bf Proof.} We will prove that all compositions in
$\widetilde{S}^+$ are trivial modulo $\widetilde{S}^+$. We consider
the following cases.

Case 1. $f=x_{mn}x_{ij}-q^{-2}x_{ij}x_{mn}, \
g=x_{ij}x_{kl}-q^{-2}x_{kl}x_{ij}, \ w=x_{mn}x_{ij}x_{kl}$.

In the case, we have
$$
(f,g)_w=-q^{-2}x_{ij}x_{mn}x_{kl}+q^{-2}x_{mn}x_{kl}x_{ij}.
$$

There are four subcases to consider.
$$
\begin{tabular}{|c|c|c|}
\hline
& $((i,j),(m,n)) \in C_1$ & $((i,j),(m,n)) \in C_3$ \\
\hline
$((k,l),(i,j)) \in C_1$ & 1.1. $((k,l),(m,n)) \in C_1$ &
1.3. $((k,l),(m,n)) \in C_4, C_5 \ or \ C_6$ \\
\hline
$((k,l),(i,j)) \in C_3$ & 1.2. $((k,l),(m,n)) \in C_4$ &
1.4. $((k,l),(m,n)) \in C_3$\\
\hline
\end{tabular}
$$

1.1. \ $((i,j),(m,n)) \in C_1, \ ((k,l),(i,j)) \in C_1$ and
$((k,l),(m,n)) \in C_1$.

Then, we have
\begin{eqnarray*}
(f,g)_w&\equiv&-q^{-4}x_{ij}x_{kl}x_{mn}+q^{-4}x_{kl}x_{mn}x_{ij}\\
&\equiv&-q^{-6}x_{kl}x_{ij}x_{mn}+q^{-6}x_{kl}x_{ij}x_{mn}\\
&\equiv&0.
\end{eqnarray*}

1.2. \ $((i,j),(m,n)) \in C_1, \ ((k,l),(i,j)) \in C_3$ and
$((k,l),(m,n)) \in C_4$.

Then, we have $(i,j)=(m,l), \ ((k,n),(i,j)) \in C_2$ and
\begin{eqnarray*}
(f,g)_w&\equiv&-q^{-2}x_{ij}[x_{kl}x_{mn}-(q^2-q^{-2})x_{kn}x_{ml}]
+q^{-2}[x_{kl}x_{mn}-(q^2-q^{-2})x_{kn}x_{ml}]x_{ij}\\
&\equiv&-q^{-4}x_{kl}x_{ij}x_{mn}+q^{-2}(q^2-q^{-2})x_{kn}x_{ij}x_{ml}\\
&&+q^{-4}x_{kl}x_{ij}x_{mn}-q^{-2}(q^2-q^{-2})x_{kn}x_{ij}x_{ml}\\
&\equiv&0.
\end{eqnarray*}

1.3. \ $((i,j),(m,n)) \in C_3, \ ((k,l),(i,j)) \in C_1$ and
$((k,l),(m,n)) \in C_4, \ C_5$ or $C_6$.

1.3.1. If $((k,l),(m,n)) \in C_4 \ (m<l)$, then $(k,n)=(i,j), \
((i,j),(m,l)) \in C_2$ and
\begin{eqnarray*}
(f,g)_w&\equiv&-q^{-2}x_{ij}[x_{kl}x_{mn}-(q^2-q^{-2})x_{kn}x_{ml}]
+q^{-2}[x_{kl}x_{mn}-(q^2-q^{-2})x_{kn}x_{ml}]x_{ij}\\
&\equiv&-q^{-4}x_{kl}x_{ij}x_{mn}+q^{-2}(q^2-q^{-2})x_{ij}x_{kn}x_{ml}
+q^{-4}x_{kl}x_{ij}x_{mn}\\
&&-q^{-2}(q^2-q^{-2})x_{kn}x_{ij}x_{ml}\\
&\equiv&0.
\end{eqnarray*}

1.3.2. If $((k,l),(m,n)) \in C_5 \ (m=l)$, then $(k,n)=(i,j)$ and
\begin{eqnarray*}
(f,g)_w&\equiv&-q^{-2}x_{ij}(q^2x_{kl}x_{mn}-qx_{kn})
+q^2(q^2x_{kl}x_{mn}-qx_{kn})x_{ij}\\
&\equiv&-x_{ij}x_{kl}x_{mn}+q^{-1}x_{ij}x_{kn}+x_{kl}x_{mn}x_{ij}
-q^{-1}x_{kn}x_{ij}\\
&\equiv&-q^{-2}x_{kl}x_{ij}x_{mn}+q^{-2}x_{kl}x_{ij}x_{mn}\\
&\equiv&0.
\end{eqnarray*}

1.3.3. If $((k,l),(m,n)) \in C_6 \ (m>l)$, then
\begin{eqnarray*}
(f,g)_w&\equiv&-q^{-2}x_{ij}x_{kl}x_{mn}+q^{-2}x_{kl}x_{mn}x_{ij}\\
&\equiv&-q^{-4}x_{kl}x_{ij}x_{mn}+q^{-4}x_{kl}x_{ij}x_{mn}\\
&\equiv&0.
\end{eqnarray*}

1.4. \ $((i,j),(m,n)) \in C_3, \ ((k,l),(i,j)) \in C_3$ and
$((k,l),(m,n)) \in C_3$.

This case is similar to 1.1.

\ \

Case 2.  $f=x_{mn}x_{ij}-q^{-2}x_{ij}x_{mn}, \
g=x_{ij}x_{kl}-x_{kl}x_{ij}, \ w=x_{mn}x_{ij}x_{kl}$.

In the case, we have
$$
(f,g)_w=-q^{-2}x_{ij}x_{mn}x_{kl}+x_{mn}x_{kl}x_{ij}.
$$

There are also four subcases to consider.
$$
\begin{tabular}{|c|c|c|}
\hline
& $((i,j),(m,n)) \in C_1$ & $((i,j),(m,n)) \in C_3$ \\
\hline
$((k,l),(i,j)) \in C_2$ & 2.1. $((k,l),(m,n)) \in C_2, C_3 \ or \ C_4$
& 2.3. $((k,l),(m,n)) \in C_2$ \\
\hline
$((k,l),(i,j)) \in C_6$ & 2.2. $((k,l),(m,n)) \in C_6$ &
2.4. $((k,l),(m,n)) \in C_6$\\
\hline
\end{tabular}
$$

2.1. \ $((i,j),(m,n)) \in C_1, \ ((k,l),(i,j)) \in C_2$ and
$((k,l),(m,n)) \in C_2, C_3$ or $C_4$.

2.1.1. If $((k,l),(m,n)) \in C_2 \ (n<l)$, then
\begin{eqnarray*}
(f,g)_w&\equiv&-q^{-2}x_{ij}x_{kl}x_{mn}+x_{kl}x_{mn}x_{ij}\\
&\equiv&-q^{-2}x_{kl}x_{ij}x_{mn}+q^{-2}x_{kl}x_{ij}x_{mn}\\
&\equiv&0.
\end{eqnarray*}

2.1.2. If $((k,l),(m,n)) \in C_3 \ (n=l)$, then
\begin{eqnarray*}
(f,g)_w&\equiv&-q^{-4}x_{ij}x_{kl}x_{mn}+q^{-2}x_{kl}x_{mn}x_{ij}\\
&\equiv&-q^{-4}x_{kl}x_{ij}x_{mn}+q^{-4}x_{kl}x_{ij}x_{mn}\\
&\equiv&0.
\end{eqnarray*}

2.1.3. If $((k,l),(m,n)) \in C_4 \ (n>l)$, then $((k,n),(i,j)) \in
C_2, \ ((i,j),(m,l)) \in C_1$ and
\begin{eqnarray*}
(f,g)_w&\equiv&-q^{-2}x_{ij}[x_{kl}x_{mn}-(q^2-q^{-2})x_{kn}x_{ml}]
+[x_{kl}x_{mn}-(q^2-q^{-2})x_{kn}x_{ml}]x_{ij}\\
&\equiv&-q^{-2}x_{kl}x_{ij}x_{mn}+q^{-2}(q^2-q^{-2})x_{kn}x_{ij}x_{ml}
+q^{-2}x_{kl}x_{ij}x_{mn}\\
&&-q^{-2}(q^2-q^{-2})x_{kn}x_{ij}x_{ml}\\
&\equiv&0.
\end{eqnarray*}

For the cases 2.2, 2.3 and 2.4, the proofs are similar to 2.1.1.

 \ \

Case 3. $f=x_{mn}x_{ij}-q^{-2}x_{ij}x_{mn}, \
g=x_{ij}x_{kl}-x_{kl}x_{ij}+(q^2-q^{-2})x_{kj}x_{il}, \
w=x_{mn}x_{ij}x_{kl}$.

In the case, we have
$$
(f,g)_w=-q^{-2}x_{ij}x_{mn}x_{kl}+x_{mn}x_{kl}x_{ij}
-(q^2-q^{-2})x_{mn}x_{kj}x_{il}.
$$

There are two subcases to consider.
$$
\begin{tabular}{|c|c|c|}
\hline
& $((i,j),(m,n)) \in C_1$& $((i,j),(m,n)) \in C_3$ \\
\hline
& 3.1. & 3.2. \\
$((k,l),(i,j)) \in C_4$ & $((k,l),(m,n)), ((k,j),(m,n)) \in C_4$ &
$((k,l),(m,n)) \in C_4, C_5 \ or \ C_6$ \\
& & $((k,j),(m,n)) \in C_3$ \\
\hline
\end{tabular}
$$

3.1. \ $((i,j),(m,n)) \in C_1, \ ((k,l),(i,j)) \in C_4$ and
$(k,l),(m,n)), ((k,j),(m,n)) \in C_4$.

Then, we have $((k,n),(i,j)) \in C_2, \ ((i,l),(m,n)) \in C_1, \
((i,l),(m,j)) \in C_1, \ ((m,l),(i,j))\\
\in C_1$ and
\begin{eqnarray*}
(f,g)_w&\equiv&-q^{-2}x_{ij}[x_{kl}x_{mn}-(q^2-q^{-2})x_{kn}x_{ml}]
+[x_{kl}x_{mn}-(q^2-q^{-2})x_{kn}x_{ml}]x_{ij}\\
&&-(q^2-q^{-2})[x_{kj}x_{mn}-(q^2-q^{-2})x_{kn}x_{mj}]x_{il}\\
&\equiv&-q^{-2}[x_{kl}x_{ij}-(q^2-q^{-2})x_{kj}x_{il}]x_{mn}
+q^{-2}(q^2-q^{-2})x_{kn}x_{ij}x_{ml}+q^{-2}x_{kl}x_{ij}x_{mn}\\
&&-(q^2-q^{-2})x_{kn}x_{ml}x_{ij}-q^{-2}(q^2-q^{-2})x_{kj}x_{il}x_{mn}
+q^{-2}(q^2-q^{-2})^2x_{kn}x_{il}x_{mj}\\
&\equiv&q^{-4}(q^2-q^{-2})x_{kn}x_{ml}x_{ij}-(q^2-q^{-2})x_{kn}x_{ml}x_{ij}
+q^{-2}(q^2-q^{-2})x_{kn}x_{ml}x_{ij}\\
&\equiv&0.
\end{eqnarray*}

3.2. \ $((i,j),(m,n)) \in C_3, \ ((k,l),(i,j)) \in C_4, \
(k,l),(m,n)) \in C_4, C_5 \ or \ C_6$ and $((k,j),(m,n)) \in C_3$.

3.2.1. If $((k,l),(m,n)) \in C_4 \ (l<m)$ and $((k,j),(m,n)) \in
C_3$, then $((k,n),(i,j)) \in C_3, \ ((i,j),(m,l)) \in C_2, \
((i,l),(m,n)) \in C_4$ and
\begin{eqnarray*}
(f,g)_w&\equiv&-q^{-2}x_{ij}[x_{kl}x_{mn}-(q^2-q^{-2})x_{kn}x_{ml}]
+[x_{kl}x_{mn}-(q^2-q^{-2})x_{kn}x_{ml}]x_{ij}\\
&&-q^{-2}(q^2-q^{-2})x_{kj}x_{mn}x_{il}\\
&\equiv&-q^{-2}[x_{kl}x_{ij}-(q^2-q^{-2})x_{kj}x_{il}]x_{mn}
+q^{-4}(q^2-q^{-2})x_{kn}x_{ij}x_{ml}+q^{-2}x_{kl}x_{ij}x_{mn}\\
&&-(q^2-q^{-2})x_{kn}x_{ij}x_{ml}-q^{-2}(q^2-q^{-2})x_{kj}[x_{il}x_{mn}
-(q^2-q^{-2})x_{in}x_{ml}]\\
&\equiv&0.
\end{eqnarray*}

3.2.2. If $((k,l),(m,n)) \in C_5 \ (l=m)$ and $((k,j),(m,n)) \in
C_3$, then $((k,l),(i,j)) \in C_4, \ ((k,n),(i,j)) \in C_3, \
((i,l),(m,n)) \in C_5$ and
\begin{eqnarray*}
(f,g)_w&\equiv&-q^{-2}x_{ij}(q^2x_{kl}x_{mn}-qx_{kn})
+(q^2x_{kl}x_{mn}-qx_{kn})x_{ij}
-q^{-2}(q^2-q^{-2})x_{kj}x_{mn}x_{il}\\
&\equiv&-[x_{kl}x_{ij}-(q^2-q^{-2})x_{kj}x_{il}]x_{mn}
+q^{-3}x_{kn}x_{ij}+x_{kl}x_{ij}x_{mn}-qx_{kn}x_{ij}\\
&&-q^{-2}(q^2-q^{-2})x_{kj}[q^2x_{il}x_{mn}-qx_{in}]\\
&\equiv&q^{-3}x_{kn}x_{ij}-qx_{kn}x_{ij}+q^{-1}(q^2-q^{-2})x_{kn}x_{ij}\\
&\equiv&0.
\end{eqnarray*}

3.2.3. If $((k,l),(m,n)) \in C_6 \ (l>m)$ and $((k,j),(m,n)) \in
C_3$, then $((i,l),(m,n)) \in C_6$ and
\begin{eqnarray*}
(f,g)_w&\equiv&-q^{-2}x_{ij}x_{kl}x_{mn}+x_{kl}x_{mn}x_{ij}
-q^{-2}(q^2-q^{-2})x_{kj}x_{mn}x_{il}\\
&\equiv&-q^{-2}[x_{kl}x_{ij}-(q^2-q^{-2})x_{kj}x_{il}]x_{mn}
+q^{-2}x_{kl}x_{ij}x_{mn}-q^{-2}(q^2-q^{-2})x_{kj}x_{il}x_{mn}\\
&\equiv&0.
\end{eqnarray*}

 \ \

Case 4. $f=x_{mn}x_{ij}-q^{-2}x_{ij}x_{mn}, \
g=x_{ij}x_{kl}-q^2x_{kl}x_{ij}+qx_{kj}, \ w=x_{mn}x_{ij}x_{kl}.$

In the case, we have
$$
(f,g)_w=-q^{-2}x_{ij}x_{mn}x_{kl}+q^2x_{mn}x_{kl}x_{ij}-qx_{mn}x_{kj}.
$$

There are two subcases to consider.
$$
\begin{tabular}{|c|c|c|}
\hline
& $((i,j),(m,n)) \in C_1$& $((i,j),(m,n)) \in C_3$ \\
\hline
& 4.1. & 4.2. \\
$((k,l),(i,j)) \in C_5$ & $((k,l),(m,n)) \in C_5$ & $((k,l),(m,n)) \in C_6$ \\
& $((k,j),(m,n)) \in C_4$ & $((k,j),(m,n)) \in C_3$ \\
\hline
\end{tabular}
$$

4.1. \ $((i,j),(m,n)) \in C_1, \ ((k,l),(i,j)) \in C_5,
((k,l),(m,n)) \in C_5$ and $((k,j),(m,n)) \in C_4$.

Then, we have $((k,n),(i,j)) \in C_2 \ (m=i)$ and
\begin{eqnarray*}
(f,g)_w&\equiv&-q^{-2}x_{ij}(q^2x_{kl}x_{mn}-qx_{kn})+q^2(q^2x_{kl}x_{mn}-qx_{kn})x_{ij}\\
&&-q[x_{kj}x_{mn}-(q^2-q^{-2})x_{kn}x_{mj}]\\
&\equiv&-(q^2x_{kl}x_{ij}-qx_{kj})x_{mn}+q^{-1}x_{kn}x_{ij}+q^2x_{kl}x_{ij}x_{mn}\\
&&-q^3x_{kn}x_{ij}-qx_{kj}x_{mn}+q(q^2-q^{-2})x_{kn}x_{mj}\\
&\equiv&0.
\end{eqnarray*}

4.2. \ $((i,j),(m,n)) \in C_3, \ ((k,l),(i,j)) \in C_5,
((k,l),(m,n)) \in C_6$ and $((k,j),(m,n)) \in C_3$.

Then, we have
\begin{eqnarray*}
(f,g)_w&\equiv&-q^{-2}x_{ij}x_{kl}x_{mn}+q^2x_{kl}x_{mn}x_{ij}-q^{-1}x_{kj}x_{mn}\\
&\equiv&-q^{-2}(q^2x_{kl}x_{ij}-qx_{kj})x_{mn}+x_{kl}x_{ij}x_{mn}-q^{-1}x_{kj}x_{mn}\\
&\equiv&0.
\end{eqnarray*}

 \ \

Case 5. $f=x_{mn}x_{ij}-x_{ij}x_{mn}, \
g=x_{ij}x_{kl}-q^{-2}x_{kl}x_{ij}, \ w=x_{mn}x_{ij}x_{kl}.$

In the case, we have
$$
(f,g)_w=-x_{ij}x_{mn}x_{kl}+q^{-2}x_{mn}x_{kl}x_{ij}.
$$

There are four subcases to consider.
$$
\begin{tabular}{|c|c|c|}
\hline
& $((i,j),(m,n)) \in C_2$ & $((i,j),(m,n)) \in C_6$ \\
\hline
$((k,l),(i,j)) \in C_1$ & 5.1. $((k,l),(m,n)) \in C_2,C_3,C_4,C_5 \ or \ C_6$ &
5.3. $((k,l),(m,n)) \in C_6$ \\
\hline
$((k,l),(i,j)) \in C_3$ & 5.2. $((k,l),(m,n)) \in C_2$ & 5.4. $((k,l),(m,n)) \in C_6$\\
\hline
\end{tabular}
$$

5.1. \ $((i,j),(m,n)) \in C_2, \ ((k,l),(i,j)) \in C_1$, and
$((k,l),(m,n)) \in C_2,C_3,C_4,C_5$ or $C_6$.

5.1.1. If $((k,l),(m,n)) \in C_2 \ (l>n)$, then we have
$((k,l),(i,j)) \in C_1$ and
\begin{eqnarray*}
(f,g)_w&\equiv&-x_{ij}x_{kl}x_{mn}+q^{-2}x_{kl}x_{mn}x_{ij}\\
&\equiv&-q^{-2}x_{kl}x_{ij}x_{mn}+q^{-2}x_{kl}x_{ij}x_{mn}\\
&\equiv&0.
\end{eqnarray*}

5.1.2. If $((k,l),(m,n)) \in C_3 \ (l=n)$, then
\begin{eqnarray*}
(f,g)_w&\equiv&-q^{-2}x_{ij}x_{kl}x_{mn}+q^{-4}x_{kl}x_{mn}x_{ij}\\
&\equiv&-q^{-4}x_{kl}x_{ij}x_{mn}+q^{-4}x_{kl}x_{ij}x_{mn}\\
&\equiv&0.
\end{eqnarray*}

5.1.3. If $((k,l),(m,n)) \in C_4 \ (m<l<n)$, then we have
$((k,l),(i,j)) \in C_1, \ ((k,n),(i,j)) \in C_1, \ ((i,j),(m,l)) \in
C_2$ and
\begin{eqnarray*}
(f,g)_w&\equiv&-x_{ij}[x_{kl}x_{mn}-(q^2-q^{-2})x_{kn}x_{ml}]
+q^{-2}[x_{kl}x_{mn}-(q^2-q^{-2})x_{kn}x_{ml}]x_{ij}\\
&\equiv&-x_{ij}x_{kl}x_{mn}+(q^2-q^{-2})x_{ij}x_{kn}x_{ml}
+q^{-2}x_{kl}x_{mn}x_{ij}-q^{-2}(q^2-q^{-2})x_{kn}x_{ml}x_{ij}\\
&\equiv&-q^{-2}x_{kl}x_{ij}x_{mn}+q^{-2}(q^2-q^{-2})x_{kn}x_{ij}x_{ml}
+q^{-2}x_{kl}x_{ij}x_{mn}\\
&&-q^{-2}(q^2-q^{-2})x_{kn}x_{ij}x_{ml}\\
&\equiv&0.
\end{eqnarray*}

5.1.4. If $((k,l),(m,n)) \in C_5 \ (m=l)$, then we have
$((k,n),(i,j)) \in C_1$ and
\begin{eqnarray*}
(f,g)_w&\equiv&-x_{ij}(q^2x_{kl}x_{mn}-qx_{kn})+q^{-2}(q^2x_{kl}x_{mn}-qx_{kn})x_{ij}\\
&\equiv&-q^2x_{ij}x_{kl}x_{mn}+qx_{ij}x_{kn}+x_{kl}x_{mn}x_{ij}-q^{-1}x_{kn}x_{ij}\\
&\equiv&-x_{kl}x_{ij}x_{mn}+q^{-1}x_{kn}x_{ij}+x_{kl}x_{ij}x_{mn}-q^{-1}x_{kn}x_{ij}\\
&\equiv&0.
\end{eqnarray*}

5.1.5. If $((k,l),(m,n)) \in C_6 \ (l<m)$, the proof is similar to
5.1.1.

For the cases of 5.2, 5.3 and 5.4, the proofs are also similar to
5.1.1.

 \ \

Case 6. $f=x_{mn}x_{ij}-x_{ij}x_{mn}, \ g=x_{ij}x_{kl}-x_{kl}x_{ij},
\ w=x_{mn}x_{ij}x_{kl}.$

In the case, we have
$$
(f,g)_w=-x_{ij}x_{mn}x_{kl}+x_{mn}x_{kl}x_{ij}.
$$

There are four subcases to consider.
$$
\begin{tabular}{|c|c|c|}
\hline
& $((i,j),(m,n)) \in C_2$ & $((i,j),(m,n)) \in C_6$ \\
\hline
$((k,l),(i,j)) \in C_2$ & 6.1. $((k,l),(m,n)) \in C_2$ &
6.3. $((k,l),(m,n)) \in C_2,C_3,C_4,C_5 \ or \ C_6$ \\
\hline
$((k,l),(i,j)) \in C_6$ & 6.2. $((k,l),(m,n)) \in C_6$ & 6.4. $((k,l),(m,n)) \in C_6$\\
\hline
\end{tabular}
$$

6.1. \ $((i,j),(m,n)) \in C_2, \ ((k,l),(m,n)) \in C_2$ and
$((k,l),(m,n)) \in C_2$.

Then, we have
\begin{eqnarray*}
(f,g)_w&\equiv&-x_{ij}x_{kl}x_{mn}+x_{kl}x_{mn}x_{ij}\\
&\equiv&-x_{kl}x_{ij}x_{mn}+x_{kl}x_{ij}x_{mn}\\
&\equiv&0.
\end{eqnarray*}

6.2. \ $((i,j),(m,n)) \in C_2, \ ((k,l),(m,n)) \in C_6$ and
$((k,l),(m,n)) \in C_6$.

This case is similar to 6.1.

6.3. \ $((i,j),(m,n)) \in C_6, \ ((k,l),(m,n)) \in C_2$ and
$((k,l),(m,n)) \in C_2,,C_3,C_4,C_5$ or $C_6$.

6.3.1. If $((k,l),(m,n)) \in C_2 \ (l>n)$, the proof is similar to
6.1.

6.3.2. If $((k,l),(m,n)) \in C_3 \ (l=n)$, then
\begin{eqnarray*}
(f,g)_w&\equiv&-q^{-2}x_{ij}x_{kl}x_{mn}+q^{-2}x_{kl}x_{mn}x_{ij}\\
&\equiv&-q^{-2}x_{kl}x_{ij}x_{mn}+q^{-2}x_{kl}x_{ij}x_{mn}\\
&\equiv&0.
\end{eqnarray*}

6.3.3. If $((k,l),(m,n)) \in C_4 \ (m<l<n)$, then we have
$((k,n),(i,j)) \in C_2, \ ((i,j),(m,n)) \in C_6$ and
\begin{eqnarray*}
(f,g)_w&\equiv&-x_{ij}[x_{kl}x_{mn}-(q^2-q^{-2})x_{kn}x_{ml}]
+[x_{kl}x_{mn}-(q^2-q^{-2})x_{kn}x_{ml}]x_{ij}\\
&\equiv&-x_{ij}x_{kl}x_{mn}+(q^2-q^{-2})x_{ij}x_{kn}x_{ml}
+x_{kl}x_{mn}x_{ij}-(q^2-q^{-2})x_{kn}x_{ml}x_{ij}\\
&\equiv&-x_{kl}x_{ij}x_{mn}+(q^2-q^{-2})x_{kn}x_{ij}x_{ml}
+x_{kl}x_{ij}x_{mn}-(q^2-q^{-2})x_{kn}x_{ij}x_{ml}\\
&\equiv&0.
\end{eqnarray*}

6.3.4. If $((k,l),(m,n)) \in C_5 \ (m=l)$, then we have
$((k,n),(i,j)) \in C_2$ and
\begin{eqnarray*}
(f,g)_w&\equiv&-x_{ij}(q^2x_{kl}x_{mn}-qx_{kn})+(q^2x_{kl}x_{mn}-qx_{kn})x_{ij}\\
&\equiv&-q^2x_{ij}x_{kl}x_{mn}+qx_{ij}x_{kn}+q^2x_{kl}x_{mn}x_{ij}-qx_{kn}x_{ij}\\
&\equiv&-q^2x_{kl}x_{ij}x_{mn}+qx_{kn}x_{ij}+q^2x_{kl}x_{ij}x_{mn}-qx_{kn}x_{ij}\\
&\equiv&0.
\end{eqnarray*}

6.3.5. If $((k,l),(m,n)) \in C_6 \ (l<m)$, the proof is similar to
6.1.

6.4. \ $((i,j),(m,n)) \in C_6, \ ((k,l),(m,n)) \in C_6$ and
$((k,l),(m,n)) \in C_6$.

This case is also similar to 6.1.

 \ \

Case 7. $f=x_{mn}x_{ij}-x_{ij}x_{mn}, \
g=x_{ij}x_{kl}-x_{kl}x_{ij}+(q^2-q^{-2})x_{kj}x_{il}, \
w=x_{mn}x_{ij}x_{kl}.$

In the case, we have
$$
(f,g)_w=-x_{ij}x_{mn}x_{kl}+x_{mn}x_{kl}x_{ij}-(q^2-q^{-2})x_{mn}x_{kj}x_{il}.
$$

There are two subcases to consider.
$$
\begin{tabular}{|c|c|c|}
\hline
 &  $((i,j),(m,n)) \in C_2$& $((i,j),(m,n)) \in C_6$ \\
\hline
 & 7.1. & 7.2. \\
$((k,l),(i,j)) \in C_4$ & $((k,l),(m,n)) \in C_2,C_3,C_4,C_5 \ or \ C_6$ &
$((k,l),(m,n)),$ \\
 & $((k,j),(m,n)) \in C_2$ & $((k,j),(m,n)) \in C_6$  \\
\hline
\end{tabular}
$$

7.1. \ $((i,j),(m,n)) \in C_2, \ ((k,l),(i,j)) \in C_4, \
((k,l),(m,n)) \in C_2,C_3,C_4,C_5$ or $C_6$ and $((k,j),(m,n)) \in
C_2$.

7.1.1. If $((k,l),(m,n)) \in C_2 \ (n<l)$ and $((k,j),(m,n)) \in
C_2$, then we have $((i,l),(m,n)) \in C_2$ and
\begin{eqnarray*}
(f,g)_w&\equiv&-x_{ij}x_{kl}x_{mn}+x_{kl}x_{mn}x_{ij}
-(q^2-q^{-2})x_{kj}x_{mn}x_{il}\\
&\equiv&-[x_{kl}x_{ij}-(q^2-q^{-2})x_{kj}x_{il}]x_{mn}
+x_{kl}x_{ij}x_{mn}-(q^2-q^{-2})x_{kj}x_{il}x_{mn}\\
&\equiv&0.
\end{eqnarray*}

7.1.2. If $((k,l),(m,n)) \in C_3 \ (n=l)$ and $((k,j),(m,n)) \in
C_2$, then  $((i,l),(m,n)) \in C_3$ and
\begin{eqnarray*}
(f,g)_w&\equiv&-q^{-2}x_{ij}x_{kl}x_{mn}+q^{-2}x_{kl}x_{mn}x_{ij}
-(q^2-q^{-2})x_{kj}x_{mn}x_{il}\\
&\equiv&-q^{-2}[x_{kl}x_{ij}-(q^2-q^{-2})x_{kj}x_{il}]x_{mn}
+q^{-2}x_{kl}x_{ij}x_{mn}-q^{-2}(q^2-q^{-2})x_{kj}x_{il}x_{mn}\\
&\equiv&0.
\end{eqnarray*}

7.1.3. If $((k,l),(m,n)) \in C_4 \ (m<l<n)$ and $((k,j),(m,n)) \in
C_2$, then we obtain $((k,l),(i,j)) \in C_4, \ ((k,n),(i,j)) \in
C_4, \ ((i,j),(m,l)) \in C_2, \ ((i,l),(m,n)) \in C_4$ and
\begin{eqnarray*}
(f,g)_w&\equiv&-x_{ij}[x_{kl}x_{mn}-(q^2-q^{-2})x_{kn}x_{ml}]
+[x_{kl}x_{mn}-(q^2-q^{-2})x_{kn}x_{ml}]x_{ij}\\
&&-(q^2-q^{-2})x_{kj}x_{mn}x_{il}\\
&\equiv&-x_{ij}x_{kl}x_{mn}+(q^2-q^{-2})x_{ij}x_{kn}x_{ml}
+x_{kl}x_{mn}x_{ij}-(q^2-q^{-2})x_{kn}x_{ml}x_{ij}\\
&&-(q^2-q^{-2})x_{kj}[x_{il}x_{mn}-(q^2-q^{-2})x_{in}x_{ml}]\\
&\equiv&-[x_{kl}x_{ij}-(q^2-q^{-2})x_{kj}x_{il}]x_{mn}
+(q^2-q^{-2})[x_{kn}x_{ij}-(q^2-q^{-2})x_{kj}x_{in}]x_{ml}\\
&&+x_{kl}x_{ij}x_{mn}-(q^2-q^{-2})x_{kn}x_{ij}x_{ml}\\
&\equiv&0.
\end{eqnarray*}

7.1.4. If $((k,l),(m,n)) \in C_5 \ (m=l)$ and $((k,j),(m,n)) \in
C_2$, then  $((k,n),(i,j)) \in C_4, \ ((i,l),(m,n)) \in C_5$ and
\begin{eqnarray*}
(f,g)_w&\equiv&-x_{ij}(q^2x_{kl}x_{mn}-qx_{kn})+(q^2x_{kl}x_{mn}
-qx_{kn})x_{ij}-(q^2-q^{-2})x_{kj}x_{mn}x_{il}\\
&\equiv&-q^2x_{ij}x_{kl}x_{mn}+qx_{ij}x_{kn}+q^2x_{kl}x_{mn}x_{ij}
-qx_{kn}x_{ij}\\
&&-(q^2-q^{-2})x_{kj}(q^2x_{il}x_{mn}-qx_{in})\\
&\equiv&-q^2[x_{kl}x_{ij}-(q^2-q^{-2})x_{kj}x_{il}]x_{mn}+q[x_{kn}x_{ij}
-(q^2-q^{-2})x_{kj}x_{in}]\\
&&+q^2x_{kl}x_{ij}x_{mn}-qx_{kn}x_{ij}-q^2(q^2-q^{-2})x_{kj}x_{il}x_{mn}
+q(q^2-q^{-2})x_{kj}x_{in}\\
&\equiv&0.
\end{eqnarray*}

7.1.5. If $((k,l),(m,n)) \in C_6 \ (l<m)$ and $((k,j),(m,n)) \in
C_2$, then $((i,l),(m,n)) \in C_6$. This case is similar to 7.1.1.

7.2. \ $((i,j),(m,n)) \in C_6, \ ((k,l),(i,j)) \in C_4, \
((k,l),(m,n)),((k,j),(m,n)) \in C_6$.

This case is also similar to 7.1.1.

 \ \

Case 8. $f=x_{mn}x_{ij}-x_{ij}x_{mn}, \
g=x_{ij}x_{kl}-q^2x_{kj}x_{il}+qx_{kj}, \ w=x_{mn}x_{ij}x_{kl}.$

In the case, we have
$$
(f,g)_w=-x_{ij}x_{mn}x_{kl}+q^2x_{mn}x_{kl}x_{ij}+qx_{mn}x_{kj}.
$$

There are two subcases to consider.
$$
\begin{tabular}{|c|c|c|}
\hline
& $((i,j),(m,n)) \in C_2$& $((i,j),(m,n)) \in C_6$ \\
\hline
& 8.1. & 8.2. \\
$((k,l),(i,j)) \in C_5$ & $((k,l),(m,n)) \in C_6$ &
$((k,l),(m,n)),((k,j),(m,n)) \in C_6$ \\
& $((k,j),(m,n)) \in C_2$ &  \\
\hline
\end{tabular}
$$

8.1. \ $((i,j),(m,n)) \in C_2, \ ((k,l),(i,j)) \in C_5, \
((k,l),(m,n)) \in C_6$ and $((k,j),(m,n)) \in C_2$. Then, we have
\begin{eqnarray*}
(f,g)_w&\equiv&-x_{ij}x_{kl}x_{mn}+q^2x_{kl}x_{mn}x_{ij}+qx_{kj}x_{mn}\\
&\equiv&-(q^2x_{kl}x_{ij}-qx_{kj})x_{mn}+q^2x_{kl}x_{ij}x_{mn}+qx_{kj}x_{mn}\\
&\equiv&0.
\end{eqnarray*}

8.2. \ $((i,j),(m,n)) \in C_6, \ ((k,l),(i,j)) \in C_5, \
((k,l),(m,n)),((k,j),(m,n)) \in C_6$.

This case is similar to 8.1.

 \ \

Case 9. $f=x_{mn}x_{ij}-x_{ij}x_{mn}+(q^2-q^{-2})x_{in}x_{mj}, \
g=x_{ij}x_{kl}-q^{-2}x_{kl}x_{ij}, \ w=x_{mn}x_{ij}x_{kl}.$

 \ \

In the case, we have
$$
(f,g)_w=-x_{ij}x_{mn}x_{kl}+(q^2-q^{-2})x_{in}x_{mj}x_{kl}+q^{-2}x_{mn}x_{kl}x_{ij}.
$$

There are two subcases to consider.
$$
\begin{tabular}{|c|c|}
\hline
& $((i,j),(m,n)) \in C_4$ \\
\hline
$((k,l),(i,j)) \in C_1$ & 9.1. $((k,l),(m,n)),((k,l),(m,j)) \in C_4, \ C_5 \ or \ C_6$ \\
\hline
$((k,l),(i,j)) \in C_3$ & 9.2. $((k,l),(m,n)) \in C_4 \ \ ((k,l),(m,j)) \in C_3$ \\
\hline
\end{tabular}
$$

9.1. \ $((i,j),(m,n)) \in C_4, \ ((k,l),(i,j)) \in C_1$ and
$((k,l),(m,n)),((k,l),(m,j)) \in C_4, \ C_5$ or $C_6$.

9.1.1. If $((k,l),(m,n)),((k,l),(m,j)) \in C_4 \ (l>m)$, then we
have $((i,j),(k,n)) \in C_1, \ ((k,n),(m,l)) \in C_2, \
((k,j),(i,n)) \in C_1, \ ((k,l),(i,n)) \in C_1, \ ((i,j),(m,l)) \in
C_2$ and
\begin{eqnarray*}
(f,g)_w&\equiv&-x_{ij}[x_{kl}x_{mn}-(q^2-q^{-2})x_{kn}x_{ml}]
+(q^2-q^{-2})x_{in}[x_{kl}x_{mj}-(q^2-q^{-2})x_{kj}x_{ml}]\\
&&+q^{-2}[x_{kl}x_{mn}-(q^2-q^{-2})x_{kn}x_{ml}]x_{ij}\\
&\equiv&-x_{ij}x_{kl}x_{mn}+(q^2-q^{-2})x_{ij}x_{kn}x_{ml}
+(q^2-q^{-2})x_{in}x_{kl}x_{mj}\\
&&-(q^2-q^{-2})^2x_{in}x_{kj}x_{ml}+q^{-2}x_{kl}x_{mn}x_{ij}
-q^{-2}(q^2-q^{-2})x_{kn}x_{ml}x_{ij}\\
&\equiv&-q^{-2}x_{kl}x_{ij}x_{mn}+(q^2-q^{-2})x_{ij}x_{kn}x_{ml}
+q^{-2}(q^2-q^{-2})x_{kl}x_{in}x_{mj}\\
&&-q^{-2}(q^2-q^{-2})^2x_{kj}x_{in}x_{ml}+q^{-2}x_{kl}[x_{ij}x_{mn}
-(q^2-q^{-2})x_{in}x_{mj}]\\
&&-q^{-2}(q^2-q^{-2})x_{kn}x_{ij}x_{ml}\\
&\equiv&(q^2-q^{-2})x_{ij}x_{kn}x_{ml}-q^{-2}(q^2-q^{-2})^2x_{kj}x_{in}x_{ml}
-q^{-4}(q^2-q^{-2})x_{ij}x_{kn}x_{ml}\\
&\equiv&0.
\end{eqnarray*}

9.1.2. If $((k,l),(m,n)),((k,l),(m,j)) \in C_5 \ (l=m)$, then we
have $((i,j),(k,n)) \in C_1, \ ((k,l),(i,n),((k,j),(i,n)) \in C_1$
and
\begin{eqnarray*}
(f,g)_w&\equiv&-x_{ij}(q^2x_{kl}x_{mn}-qx_{kn})+(q^2-q^{-2})x_{in}(q^2x_{kl}x_{mj}-qx_{kj})\\
&&+q^{-2}(q^2x_{kl}x_{mn}-qx_{kn})x_{ij}\\
&\equiv&-q^2x_{ij}x_{kl}x_{mn}+qx_{ij}x_{kn}
+q^2(q^2-q^{-2})x_{in}x_{kl}x_{mj}-q(q^2-q^{-2})x_{in}x_{kj}\\
&&+x_{kl}x_{mn}x_{ij}-q^{-1}x_{kn}x_{ij}\\
&\equiv&-x_{kl}x_{ij}x_{mn}+qx_{ij}x_{kn}
+(q^2-q^{-2})x_{kl}x_{in}x_{mj}-q^{-1}(q^2-q^{-2})x_{kj}x_{in}\\
&&+x_{kl}[x_{ij}x_{mn}-(q^2-q^{-2})x_{in}x_{mj}]-q^{-3}x_{ij}x_{kn}\\
&\equiv&qx_{ij}x_{kn}-qx_{kj}x_{in}+q^{-3}x_{kj}x_{in}-q^{-3}x_{ij}x_{kn}\\
&\equiv&0.
\end{eqnarray*}

9.1.3. If $((k,l),(m,n)),((k,l),(m,j)) \in C_6 \ (l<m)$, then we
have $((k,l),(i,n)) \in C_1$ and
\begin{eqnarray*}
(f,g)_w&\equiv&-x_{ij}x_{kl}x_{mn}-(q^2-q^{-2})x_{in}x_{kl}x_{mj}+q^{-2}x_{kl}x_{mn}x_{ij}\\
&\equiv&-q^{-2}x_{kl}x_{ij}x_{mn}-q^{-2}(q^2-q^{-2})x_{kl}x_{in}x_{mj}
+q^{-2}x_{kl}[x_{ij}x_{mn}-(q^2-q^{-2})x_{in}x_{mj}]\\
&\equiv&0.
\end{eqnarray*}

9.2. \ $((i,j),(m,n)) \in C_4, \ ((k,l),(i,j)) \in C_3$,
$((k,l),(m,n))  \in C_4$ and $((k,l),(m,j)) \in C_3$. Then, we have
$((k,n),(i,j)) \in C_2, \ ((k,l),(i,n)) \in C_4, \ ((i,j),(m,l)) \in
C_3$ and
\begin{eqnarray*}
(f,g)_w&\equiv&-x_{ij}[x_{kl}x_{mn}-(q^2-q^{-2})x_{kn}x_{ml}]
+q^{-2}(q^2-q^{-2})x_{in}x_{kl}x_{mj}\\
&&+q^{-2}[x_{kl}x_{mn}-(q^2-q^{-2})x_{kn}x_{ml}]x_{ij}\\
&\equiv&-x_{ij}x_{kl}x_{mn}+(q^2-q^{-2})x_{ij}x_{kn}x_{ml}
+q^{-2}(q^2-q^{-2})x_{in}x_{kl}x_{mj}+q^{-2}x_{kl}x_{mn}x_{ij}\\
&&-q^{-2}(q^2-q^{-2})x_{kn}x_{ml}x_{ij}\\
&\equiv&-q^{-2}x_{kl}x_{ij}x_{mn}+(q^2-q^{-2})x_{kn}x_{ij}x_{ml}
+q^{-2}(q^2-q^{-2})[x_{kl}x_{in}\\
&&-(q^2-q^{-2})x_{kn}x_{il}]x_{mj}+q^{-2}x_{kl}[x_{ij}x_{mn}-(q^2-q^{-2})x_{in}x_{mj}]\\
&&-q^{-4}(q^2-q^{-2})x_{kn}x_{ij}x_{ml}\\
&\equiv&(q^2-q^{-2})x_{kn}x_{ij}x_{ml}-q^{-2}(q^2-q^{-2})x_{kn}x_{il}x_{mj}
-q^{-4}(q^2-q^{-2})x_{kn}x_{ij}x_{ml}\\
&\equiv&0.
\end{eqnarray*}

 \ \

Case 10. $f=x_{mn}x_{ij}-x_{ij}x_{mn}+(q^2-q^{-2})x_{in}x_{mj}, \
g=x_{ij}x_{kl}-x_{kl}x_{ij}, \ w=x_{mn}x_{ij}x_{kl}.$

In the case, we have
$$
(f,g)_w=-x_{ij}x_{mn}x_{kl}+(q^2-q^{-2})x_{in}x_{mj}x_{kl}+x_{mn}x_{kl}x_{ij}.
$$

There are two subcases to consider.
$$
\begin{tabular}{|c|c|}
\hline
& $((i,j),(m,n)) \in C_4$ \\
\hline
$((k,l),(i,j)) \in C_2$ &
10.1. $((k,l),(m,n)) \in C_2,C_3 \ or \ C_4 \ \ ((k,l),(m,j)) \in C_2$ \\
\hline
$((k,l),(i,j)) \in C_6$ & 10.2. $((k,l),(m,n)),((k,l),(m,j)) \in C_6$ \\
\hline
\end{tabular}
$$

10.1. \ $((i,j),(m,n)) \in C_4, \ ((k,l),(i,j)) \in C_2, \
((k,l),(m,n)) \in C_2,C_3$ or $C_4$ and $((k,l),(m,j)) \in C_2$.

10.1.1. If $((k,l),(m,n)) \in C_2 \ (l>n)$, then we have
$((k,l),(i,n)) \in C_2$ and
\begin{eqnarray*}
(f,g)_w&\equiv&-x_{ij}x_{kl}x_{mn}+(q^2-q^{-2})x_{in}x_{kl}x_{mj}+x_{kl}x_{mn}x_{ij}\\
&\equiv&-x_{kl}x_{ij}x_{mn}+(q^2-q^{-2})x_{kl}x_{in}x_{mj}+x_{kl}[x_{ij}x_{mn}
-(q^2-q^{-2})x_{in}x_{mj}]\\
&\equiv&0.
\end{eqnarray*}

10.1.2. If $((k,l),(m,n)) \in C_3 \ (l=n)$, then we have
$((k,l),(i,n)) \in C_3$ and
\begin{eqnarray*}
(f,g)_w&\equiv&-q^{-2}x_{ij}x_{kl}x_{mn}+(q^2-q^{-2})x_{in}x_{kl}x_{mj}
+q^{-2}x_{kl}x_{mn}x_{ij}\\
&\equiv&-q^{-2}x_{kl}x_{ij}x_{mn}+q^{-2}(q^2-q^{-2})x_{kl}x_{in}x_{mj}
+q^{-2}x_{kl}[x_{ij}x_{mn}-(q^2-q^{-2})x_{in}x_{mj}]\\
&\equiv&0.
\end{eqnarray*}

10.1.3. If $((k,l),(m,n)) \in C_4 \ (l<n)$, then we have
$((k,n),(i,j)) \in C_2, \ ((k,l),(i,n)) \in C_4, \ ((i,j),(m,l)) \in
C_4$ and
\begin{eqnarray*}
(f,g)_w&\equiv&-x_{ij}[x_{kl}x_{mn}-(q^2-q^{-2})x_{kn}x_{ml}]
+(q^2-q^{-2})x_{in}x_{kl}x_{mj}\\
&&+[x_{kl}x_{mn}-(q^2-q^{-2})x_{kn}x_{ml}]x_{ij}\\
&\equiv&-x_{ij}x_{kl}x_{mn}+(q^2-q^{-2})x_{ij}x_{kn}x_{ml}
+(q^2-q^{-2})x_{in}x_{kl}x_{mj}\\
&&+x_{kl}x_{mn}x_{ij}-(q^2-q^{-2})x_{kn}x_{ml}x_{ij}\\
&\equiv&-x_{kl}x_{ij}x_{mn}+(q^2-q^{-2})x_{kn}x_{ij}x_{ml}
+(q^2-q^{-2})[x_{kl}x_{in}-(q^2-q^{-2})x_{kn}x_{il}]x_{mj}\\
&&+x_{kl}[x_{ij}x_{mn}-(q^2-q^{-2})x_{in}x_{mj}]
-(q^2-q^{-2})x_{kn}[x_{ij}x_{ml}-(q^2-q^{-2})x_{il}x_{mj}]\\
&\equiv&0.
\end{eqnarray*}

10.2. \ $((i,j),(m,n)) \in C_4, \ ((k,l),(i,j)) \in C_6, \
((k,l),(m,n)), (k,l),(m,j)) \in C_6$.

This case is similar to 10.1.

 \ \

Case 11. $f=x_{mn}x_{ij}-x_{ij}x_{mn}+(q^2-q^{-2})x_{in}x_{mj}, \
g=x_{ij}x_{kl}-x_{kl}x_{ij}+(q^2-q^{-2})x_{kj}x_{il}, \
w=x_{mn}x_{ij}x_{kl}.$

In the case, we have
$$
(f,g)_w=-x_{ij}x_{mn}x_{kl}+(q^2-q^{-2})x_{in}x_{mj}x_{kl}+x_{mn}x_{kl}x_{ij}
-(q^2-q^{-2})x_{mn}x_{kj}x_{il},
$$
with
$$
\begin{tabular}{|c|c|}
\hline
& $((i,j),(m,n)) \in C_4$ \\
\hline
$((k,l),(i,j)) \in C_4$ & $((k,l),(m,n)),((k,l),(m,j)) \in C_4,C_5 \ or \ C_6$ \\
\hline
\end{tabular}
$$

11.1. If $((k,l),(m,n)), \ ((k,l),(m,j)) \in C_4 \ \ (l>m)$, then we
have $((k,n),(i,j)) \in C_2, \ ((k,l),(i,n)) \in C_4, \
((k,j),(i,n)) \in C_4, \ ((i,j),(m,l)) \in C_2, \ ((i,l),(m,n)) \in
C_4,\\ ((i,l),(m,j)) \in C_4$ and
\begin{eqnarray*}
(f,g)_w&\equiv&-x_{ij}[x_{kl}x_{mn}-(q^2-q^{-2})x_{kn}x_{ml}]
+(q^2-q^{-2})x_{in}[x_{kl}x_{mj}-(q^2-q^{-2})x_{kj}x_{ml}]\\
&&+[x_{kl}x_{mn}-(q^2-q^{-2})x_{kn}x_{ml}]x_{ij}
-(q^2-q^{-2})[x_{kj}x_{mn}-(q^2-q^{-2})x_{kn}x_{mj}]x_{il}\\
&\equiv&-x_{ij}x_{kl}x_{mn}+(q^2-q^{-2})x_{ij}x_{kn}x_{ml}
+(q^2-q^{-2})x_{in}x_{kl}x_{mj}\\
&&-(q^2-q^{-2})x_{in}x_{kj}x_{ml}+x_{kl}x_{mn}x_{ij}-(q^2-q^{-2})x_{kn}x_{ml}x_{ij}\\
&&-(q^2-q^{-2})x_{kj}x_{mn}x_{il}+(q^2-q^{-2})^2x_{kn}x_{mj}x_{il}\\
&\equiv&-[x_{kl}x_{ij}-(q^2-q^{-2})x_{kj}x_{il}]x_{mn}
+(q^2-q^{-2})x_{kn}x_{ij}x_{ml}+(q^2-q^{-2})[x_{kj}x_{in}\\
&&-(q^2-q^{-2})x_{kn}x_{il}]x_{mj}-(q^2-q^{-2})[x_{kj}x_{in}
-(q^2-q^{-2})x_{kn}x_{ij}]x_{ml}\\
&&+x_{kl}[x_{ij}x_{mn}-(q^2-q^{-2})x_{in}x_{mj}]-(q^2-q^{-2})x_{kn}x_{ij}x_{ml} \\
&&-(q^2-q^{-2})x_{kj}[x_{il}x_{mn}-(q^2-q^{-2})x_{in}x_{ml}]\\
&&+(q^2-q^{-2})^2x_{kn}[x_{il}x_{mj}-(q^2-q^{-2})x_{ij}x_{ml}]\\
&\equiv&0.
\end{eqnarray*}

11.2. If $((k,l),(m,n)), \ ((k,l),(m,j)) \in C_5 \ (l=m)$, then we
have $((k,n),(i,j)) \in C_2, \ ((k,l),(i,n)) \in C_4, \
((k,j),(i,n)) \in C_4, \ ((i,l),(m,n)) \in C_5, \ ((i,l),(m,j)) \in
C_5$ and
\begin{eqnarray*}
(f,g)_w&\equiv&-x_{ij}(q^2x_{kl}x_{mn}-qx_{kn})+(q^2-q^{-2})x_{in}(q^2x_{kl}x_{mj}
-qx_{kj})+(q^2x_{kl}x_{mn}-qx_{kn})x_{ij}\\
&&-(q^2-q^{-2})[x_{kj}x_{mn}-(q^2-q^{-2})x_{kn}x_{mj}]x_{il}\\
&\equiv&-q^2x_{ij}x_{kl}x_{mn}+qx_{ij}x_{kn}
+q^2(q^2-q^{-2})x_{in}x_{kl}x_{mj}-q(q^2-q^{-2})x_{in}x_{kj}\\
&&+q^2x_{kl}x_{mn}x_{ij}-qx_{kn}x_{ij}-(q^2-q^{-2})x_{kj}x_{mn}x_{il}
+(q^2-q^{-2})^2x_{kn}x_{mj}x_{il}\\
&\equiv&-q^2[x_{kl}x_{ij}-(q^2-q^{-2})x_{kj}x_{il}]x_{mn}+qx_{kn}x_{ij}
+q^2(q^2-q^{-2})[x_{kl}x_{in}\\
&&-(q^2-q^{-2})x_{kn}x_{il}]x_{mj}-q(q^2-q^{-2})[x_{kj}x_{in}-(q^2-q^{-2})x_{kn}x_{ij}]\\
&&+q^2x_{kl}[x_{ij}x_{mn}-(q^2-q^{-2})x_{in}x_{mj}]-qx_{kn}x_{ij}\\
&&-(q^2-q^{-2})x_{kj}[q^2x_{il}x_{mn}-qx_{in}]
+(q^2-q^{-2})^2x_{kn}[q^2x_{il}x_{mj}-qx_{ij}]\\
&\equiv&0.
\end{eqnarray*}

11.3.  If $((k,l),(m,n)), \ ((k,l),(m,j)) \in C_6 \ (l<m)$, then
$((k,j),(m,n)) \in C_4, \ ((k,l),(i,n)) \in C_4, \
((i,l),(m,n)),((i,l),(m,j)) \in C_6$ and
\begin{eqnarray*}
(f,g)_w&\equiv&-x_{ij}x_{kl}x_{mn}+(q^2-q^{-2})x_{in}x_{kl}x_{mj}+x_{kl}x_{mn}x_{ij}\\
&&-(q^2-q^{-2})[x_{kj}x_{mn}-(q^2-q^{-2})x_{kn}x_{mj}]x_{il}\\
&\equiv&-[x_{kl}x_{ij}-(q^2-q^{-2})x_{kj}x_{il}]x_{mn}
+(q^2-q^{-2})[x_{kl}x_{in}-(q^2-q^{-2})x_{kn}x_{il}]x_{mj}\\
&&+x_{kl}[x_{ij}x_{mn}-(q^2-q^{-2})x_{in}x_{mj}]-(q^2-q^{-2})x_{kj}x_{il}x_{mn}
+(q^2-q^{-2})^2x_{kn}x_{il}x_{mj}\\
&\equiv&0.
\end{eqnarray*}

 \ \

Case 12. $f=x_{mn}x_{ij}-x_{ij}x_{mn}+(q^2-q^{-2})x_{in}x_{mj}, \
g=x_{ij}x_{kl}-q^2x_{kl}x_{ij}+qx_{kj}, \ w=x_{mn}x_{ij}x_{kl},$
with
$$
\begin{tabular}{|c|c|}
\hline
& $((i,j),(m,n)) \in C_4$ \\
\hline
$((k,l),(i,j)) \in C_5$ & $((k,l),(m,n)),((k,l),(m,j)) \in C_6$ \\
& $((k,j),(m,n)) \in C_4 \ \ ((k,l),(i,n)) \in C_5$\\
\hline
\end{tabular}
$$
In the case, we can deduce that
\begin{eqnarray*}
(f,g)_w&=&-x_{ij}x_{mn}x_{kl}+(q^2-q^{-2})x_{in}x_{mj}x_{kl}
+q^2x_{mn}x_{kl}x_{ij}-qx_{mn}x_{kj}\\
&\equiv&-x_{ij}x_{kl}x_{mn}+(q^2-q^{-2})x_{in}x_{kl}x_{mj}+q^2x_{kj}x_{mn}x_{ij}\\
&&-q[x_{kj}x_{mn}-(q^2-q^{-2})x_{kn}x_{mj}]\\
&\equiv&-(q^2x_{kl}x_{ij}-qx_{kj})x_{mn}+(q^2-q^{-2})(q^2x_{kl}x_{in}-qx_{kn})x_{mj}\\
&&+q^2x_{kl}[x_{ij}x_{mn}-(q^2-q^{-2})x_{in}x_{mj}]-qx_{kj}x_{mn}
+q(q^2-q^{-2})x_{kn}x_{mj}\\
&\equiv&-q^2x_{kl}x_{ij}x_{mn}+qx_{kj}x_{mn}+q^2(q^2-q^{-2})x_{kl}x_{in}x_{mj}
-q(q^2-q^{-2})x_{kn}x_{mj}\\
&&+q^2x_{kl}x_{ij}x_{mn}-q^2(q^2-q^{-2})x_{kl}x_{in}x_{mj}-qx_{kj}x_{mn}
+q(q^2-q^{-2})x_{kn}x_{mj}\\
&\equiv&0.
\end{eqnarray*}

\ \

Case 13. $f=x_{mn}x_{ij}-q^2x_{ij}x_{mn}+qx_{in}, \
g=x_{ij}x_{kl}-q^{-2}x_{kl}x_{ij}, \ w=x_{mn}x_{ij}x_{kl}.$

In the case, we have
$$
(f,g)_w=-q^2x_{ij}x_{mn}x_{kl}+qx_{in}x_{kl}+q^{-2}x_{mn}x_{kl}x_{ij}.
$$

There are two subcases to consider.
$$
\begin{tabular}{|c|c|}
\hline
& $((i,j),(m,n)) \in C_5$ \\
\hline
$((k,l),(i,j)) \in C_1$ & 13.1. $((k,l),(m,n)) \in C_6 \ \ ((k,l),(i,n)) \in C_1$ \\
\hline
$((k,l),(i,j)) \in C_3$ & 13.2. $((k,l),(m,n)) \in C_5 \ \ ((k,l),(i,n)) \in C_4$ \\
\hline
\end{tabular}
$$

13.1. \ $((i,j),(m,n)) \in C_5, \ ((k,l),(i,j)) \in C_1, \
((k,l),(m,n)) \in C_6$ and $((k,l),(i,n)) \in C_1$. Then, we have
\begin{eqnarray*}
(f,g)_w&=&-q^2x_{ij}x_{kl}x_{mn}+q^{-1}x_{kl}x_{in}+q^{-2}x_{kl}x_{mn}x_{ij}\\
&\equiv&-x_{kl}x_{ij}x_{mn}+q^{-1}x_{kl}x_{in}+q^{-2}x_{kl}(q^2x_{ij}x_{mn}-qx_{in})\\
&\equiv&-x_{kl}x_{ij}x_{mn}+q^{-1}x_{kl}x_{in}+x_{kl}x_{ij}x_{mn}-q^{-1}x_{kl}x_{in}\\
&\equiv&0.
\end{eqnarray*}

13.2. \ $((i,j),(m,n)) \in C_5, \ ((k,l),(i,j)) \in C_3, \
((k,l),(m,n)) \in C_5$ and $((k,l),(i,n)) \in C_4$. Then, we have
$((k,n),(i,j)) \in C_2$ and
\begin{eqnarray*}
(f,g)_w&\equiv&-q^2x_{ij}(q^2x_{kl}x_{mn}-qx_{kn})+q[x_{kl}x_{in}
-(q^2-q^{-2})x_{kn}x_{il}]\\
&&-q^{-2}(q^2x_{kl}x_{mn}-qx_{kn})x_{ij}\\
&\equiv&-q^4x_{ij}x_{kl}x_{mn}+q^3x_{ij}x_{kn}+qx_{kl}x_{in}
-q(q^2-q^{-2})x_{kn}x_{il}+x_{kl}x_{mn}x_{ij}\\
&&-q^{-1}x_{kn}x_{ij}\\
&\equiv&-q^2x_{kl}x_{ij}x_{mn}+q^3x_{kn}x_{ij}+qx_{kl}x_{in}-q^3x_{kn}x_{il}\\
&&+q^{-1}x_{kn}x_{il}+q^2x_{kl}x_{ij}x_{mn}-qx_{kl}x_{in}-q^{-1}x_{kn}x_{ij}\\
&\equiv&0.
\end{eqnarray*}

\ \

Case 14. $f=x_{mn}x_{ij}-q^2x_{ij}x_{mn}+qx_{in}, \
g=x_{ij}x_{kl}-x_{kl}x_{ij}, \ w=x_{mn}x_{ij}x_{kl}.$

In the case, we have
$$
(f,g)_w=-q^2x_{ij}x_{mn}x_{kl}+qx_{in}x_{kl}+x_{mn}x_{kl}x_{ij}.
$$

There are two subcases to consider.
$$
\begin{tabular}{|c|c|}
\hline
& $((i,j),(m,n)) \in C_5$ \\
\hline
$((k,l),(i,j)) \in C_2$ & 14.1. $((k,l),(m,n)),((k,l),(i,n)) \in C_2,C_3 \ or \ C_4$ \\
\hline
$((k,l),(i,j)) \in C_6$ & 14.2. $((k,l),(m,n)),((k,l),(i,n)) \in C_6$ \\
\hline
\end{tabular}
$$

14.1. \ $((i,j),(m,n)) \in C_5, \ ((k,l),(i,j)) \in C_2$ and
$((k,l),(m,n)),((k,l),(i,n)) \in C_2,C_3$ or $C_4$.

14.1.1. If $((k,l),(m,n))$ and $((k,l),(i,n)) \in C_2 \ (l>n)$, then
\begin{eqnarray*}
(f,g)_w&=&-q^2x_{ij}x_{kl}x_{mn}+qx_{kl}x_{in}+x_{kl}x_{mn}x_{ij}\\
&\equiv&-q^2x_{kl}x_{ij}x_{mn}+qx_{kl}x_{in}+x_{kl}(q^2x_{ij}x_{mn}-qx_{in})\\
&\equiv&0.
\end{eqnarray*}

14.1.2. If $((k,l),(m,n))$ and $((k,l),(i,n)) \in C_3 \ (l=n)$, then
\begin{eqnarray*}
(f,g)_w&=&-x_{ij}x_{kl}x_{mn}+q^{-1}x_{kl}x_{in}+q^{-2}x_{kl}x_{mn}x_{ij}\\
&\equiv&-x_{kl}x_{ij}x_{mn}+q^{-1}x_{kl}x_{in}+x_{kl}x_{ij}x_{mn}-q^{-1}x_{kl}x_{in}\\
&\equiv&0.
\end{eqnarray*}

14.1.3. If $((k,l),(m,n)), \ ((k,l),(i,n)) \in C_4 \ (l<n)$, then we
have $((k,n),(i,j)) \in C_2, \ ((i,j),(m,l)) \in C_5$ and
\begin{eqnarray*}
(f,g)_w&\equiv&-q^2x_{ij}[x_{kl}x_{mn}-(q^2-q^{-2})x_{kn}x_{ml}]
+q[x_{kl}x_{in}-(q^2-q^{-2})x_{kn}x_{il}]\\
&&+[x_{kl}x_{mn}-(q^2-q^{-2})x_{kn}x_{ml}]x_{ij}\\
&\equiv&-q^2x_{ij}x_{kl}x_{mn}+q^2(q^2-q^{-2})x_{ij}x_{kn}x_{ml}
+qx_{kl}x_{in}-q(q^2-q^{-2})x_{kn}x_{il}\\
&&+x_{kl}x_{mn}x_{ij}-(q^2-q^{-2})x_{kn}x_{ml}x_{ij}\\
&\equiv&-q^2x_{kl}x_{ij}x_{mn}+q^2(q^2-q^{-2})x_{kn}x_{ij}x_{ml}
+qx_{kl}x_{in}-q(q^2-q^{-2})x_{kn}x_{il}\\
&&+x_{kn}(q^2x_{ij}x_{mn}-qx_{in})-(q^2-q^{-2})x_{kn}(q^2x_{ij}x_{mn}-qx_{il})\\
&\equiv&0.
\end{eqnarray*}

14.2. \ $((i,j),(m,n)) \in C_5, \ ((k,l),(i,j)) \in C_6$ and
$((k,l),(m,n)),((k,l),(i,n)) \in C_6$.

This case is similar to 14.1.1.

 \ \

Case 15. $f=x_{mn}x_{ij}-q^2x_{ij}x_{mn}+qx_{in}, \
g=x_{ij}x_{kl}-x_{kl}x_{ij}+(q^2-q^{-2})x_{kj}x_{il}, \
w=x_{mn}x_{ij}x_{kl},$ with
$$
\begin{tabular}{|c|c|}
\hline
& $((i,j),(m,n)) \in C_5$ \\
\hline
$((k,l),(i,j)) \in C_4$ & $((k,l),(m,n)) \in C_6 \ \ ((k,l),(i,n)) \in C_4$ \\
& $((k,j),(m,n)) \in C_5 \ \ ((i,l),(m,n)) \in C_6$\\
\hline
\end{tabular}
$$
Then, we have
\begin{eqnarray*}
(f,g)_w&=&-q^2x_{ij}x_{mn}x_{kl}+qx_{in}x_{kl}+x_{mn}x_{kl}x_{ij}
-(q^2-q^{-2})x_{mn}x_{kj}x_{il}\\
&\equiv&-q^2x_{ij}x_{kl}x_{mn}+q[x_{kl}x_{in}
-(q^2-q^{-2})x_{kn}x_{il}]+x_{kl}x_{mn}x_{ij}\\
&&-(q^2-q^{-2})(q^2x_{kj}x_{mn}-qx_{kn})x_{il}\\
&\equiv&-q^2[x_{kl}x_{ij}-(q^2-q^{-2})x_{kj}x_{il}]x_{mn}
+qx_{kl}x_{in}-q(q^2-q^{-2})x_{kn}x_{il}\\
&&+x_{kl}(q^2x_{ij}x_{mn}-qx_{in})-q^2(q^2-q^{-2})x_{kj}x_{mn}x_{il}
+q(q^2-q^{-2})x_{kn}x_{il}\\
&\equiv&-q^2x_{kl}x_{ij}x_{mn}+q^2(q^2-q^{-2})x_{kj}x_{il}x_{mn}
+qx_{kl}x_{in}-q(q^2-q^{-2})x_{kn}x_{il}\\
&&+q^2x_{kl}x_{ij}x_{mn}-qx_{kl}x_{in}-q^2(q^2-q^{-2})x_{kj}x_{il}x_{mn}
+q(q^2-q^{-2})x_{kn}x_{il}\\
&\equiv&0.
\end{eqnarray*}

\ \

Case 16. $f=x_{mn}x_{ij}-q^2x_{ij}x_{mn}+qx_{in}, \
g=x_{ij}x_{kl}-q^2x_{kl}x_{ij}+qx_{kj}, \ w=x_{mn}x_{ij}x_{kl},$
with
$$
\begin{tabular}{|c|c|}
\hline
& $((i,j),(m,n)) \in C_5$ \\
\hline
$((k,l),(i,j)) \in C_5$ &
$((k,l),(m,n)) \in C_6 \ \ ((k,l),(i,n)),((k,j),(m,n)) \in C_5$ \\
\hline
\end{tabular}
$$
In the case, we have
\begin{eqnarray*}
(f,g)_w&=&-q^2x_{ij}x_{mn}x_{kl}+qx_{in}x_{kl}+q^2x_{mn}x_{kl}x_{ij}-qx_{mn}x_{kj}\\
&\equiv&-q^2x_{ij}x_{kl}x_{mn}+q(q^2x_{kl}x_{in}-qx_{kn})
+q^2x_{kl}x_{mn}x_{ij}-q(q^2x_{kj}x_{mn}-qx_{kn})\\
&\equiv&-q^2(q^2x_{kl}x_{ij}-qx_{kj})x_{mn}+q^3x_{kl}x_{in}
-q^2x_{kn}+q^2x_{kl}(q^2x_{ij}x_{mn}-qx_{in})\\
&&-q^3x_{kj}x_{mn}+q^2x_{kn}\\
&\equiv&0.
\end{eqnarray*}

Thus, $\widetilde{S}^+$ is a Gr\"{o}bner-Shirshov basis. This
completes the proof of Theorem \ref{t1}. \ \ $\square$

\ \

Similarly, with the deg-lex order on $\widetilde{Y}^*$,
$\widetilde{S}^-$ is a Gr\"{o}bner-Shirshov basis for
$U^-_q(A_N)=k\langle \widetilde{Y}|\widetilde{S}^-\rangle$.

We now use the same notation as before. Order the generators by:
$x_i>x_j, \ h_i>h^{-1}_i>h_j>h^{-1}_j, \ y_i>y_j$ if $i>j$, and
$x_i>h_j^{\pm1}>y_m$ for all $i, \ j, \ m$. Then we obtain a well
order (deg-lex) on $\widetilde{X}\cup H\cup \widetilde{Y}$. Thus, by
Theorem \ref{t1}, we re-obtain the following theorem in \cite{BoMa}.

\begin{theorem} (\cite{BoMa} Theorem 2.7)
Let the notation be as before. Then with the deg-lex order on
$\{\widetilde{X}\cup H\cup \widetilde{Y}\}^*$, $\widetilde{S}^+\cup
T \cup K \cup \widetilde{S}^-$ is a Gr\"{o}bner-Shirshov basis for
$U_q(A_N)=k \langle \widetilde{X}\cup H\cup
\widetilde{Y}|\widetilde{S}^+\cup T \cup K \cup
\widetilde{S}^-\rangle$.
\end{theorem}

\ \

\noindent{\bf Acknowledgement}: The authors would like to
express their deepest gratitude to Professor L. A. Bokut for his
kind guidance, useful discussions and enthusiastic encouragement.

\end{document}